\documentclass[a4paper,12pt,twoside,leqno,final]{amsart}
\usepackage{amsmath}
\usepackage{amssymb}

\newtheorem{thm}{Theorem}[section]
\newtheorem{lem}[thm]{Lemma}
\newtheorem{defn}[thm]{Definition}

\numberwithin{equation}{section}

\title[Univalent polynomials and Koebe's 1/4 theorem]
{Univalent polynomials\\ 
and Koebe's one-quarter theorem}
\author{Dmitriy Dmitrishin}
\address{Odessa National Polytechnic University, 1 Shevchenko Ave., Odessa 65044, Ukraine}
\email{dmitrishin@opu.ua}
\author{Konstantin Dyakonov}
\address{Departament de Matem\`atiques i Inform\`atica, IMUB, BGSMath, Universitat de Barcelona, Gran Via 585, E-08007 Barcelona, Spain}
\address{ICREA, Pg. Llu\'is Companys 23, E-08010 Barcelona, Spain}
\email{konstantin.dyakonov@icrea.cat}
\author{Alex Stokolos}
\address{Department of Mathematical Sciences, Georgia Southern University, Statesboro, GA 30460, USA}
\email{astokolos@georgiasouthern.edu}
\thanks{The second author was supported in part by grants MTM2014-51834-P and MTM2017-83499-P from El Ministerio de Econom\'ia y Competitividad (Spain), and by grant 2017-SGR-358 from AGAUR (Generalitat de Catalunya).}

\begin{document}
\begin{abstract}
The famous Koebe $\frac14$ theorem deals with univalent (i.e., injective) analytic functions $f$ on the unit disk $\mathbb D$. It states that if $f$ is normalized so that $f(0)=0$ and $f'(0)=1$, then the image $f(\mathbb D)$ contains the disk of radius $\frac14$ about the origin, the value $\frac14$ being best possible. Now suppose $f$ is only allowed to range over the univalent polynomials of some fixed degree. What is the optimal radius in the Koebe-type theorem that arises? And for which polynomials is it attained? A plausible conjecture is stated, and the case of small degrees is settled.
\end{abstract}

\maketitle

\section{Introduction}

Suppose you can solve a certain problem that involves general analytic functions, perhaps lying in some (fairly large) class. Then you restrict your attention to the set of polynomials of a fixed degree that are in the same class. Can you also solve the restricted (polynomial) version of the problem that arises? Well, not necessarily. While many a problem is sure to simplify or trivialize completely, there are others that become dramatically harder. In what follows, we deal with a situation of the latter kind. 

\par Our starting point is the classical Koebe one-quarter theorem, a cornerstone of geometric function theory. Recall, first of all, that an analytic function on a domain $\Omega\subset\mathbb C$ is said to be {\it univalent} if it is one-to-one, i.e., takes distinct values at distinct points of $\Omega$. Now let $S$ denote the set of univalent functions $f$ on the disk $\mathbb D:=\{z\in\mathbb C:|z|<1\}$ that have a Taylor series expansion of the form 
\begin{equation}\label{eqn:taylor}
f(z)=z+\sum_{n=2}^\infty a_nz^n
\end{equation}
(so that $f(0)=0$ and $f'(0)=1$). The one-quarter theorem -- which was actually conjectured by Koebe in 1907 and proved somewhat later by Bieberbach -- reads as follows. 

\begin{thm}[Koebe's $\frac14$ theorem]\label{thm:quarter}
For every $f\in S$, the range $f(\mathbb D)$ contains the disk $\{w:|w|<\frac14\}$.
\end{thm}

See, e.g., \cite[Chapter 2]{Du} or \cite[Chapter 14]{R} for a proof. Furthermore, since the so-called {\it Koebe function} 
\begin{equation}\label{eqn:koebefun}
\mathcal K(z):=\frac z{(1-z)^2}=z+2z^2+3z^3+\dots,\qquad z\in\mathbb D,
\end{equation}
belongs to $S$ and maps $\mathbb D$ onto the slit plane $\mathbb C\setminus(-\infty,-\frac14]$, we see that the radius $\frac14$ in Koebe's theorem is optimal; indeed, no larger number would do. 

\par Motivated by this result, we introduce the following notation. Given a set $X\subset S$, we write $\rho(X)$ for the supremum of those $r>0$ for which the common range $\bigcap_{f\in X}f(\mathbb D)$ contains the disk $\{w:|w|<r\}$. The number $\rho(X)$ will be referred to as the {\it Koebe radius} for $X$. Clearly, we always have $\rho(X)\ge\frac14$. 

\par An important subclass of $S$, to be denoted by $S_{\mathbb R}$, is the set of univalent functions $f$ of the form \eqref{eqn:taylor} whose coefficients $a_n$ are all real. Because the Koebe function $\mathcal K$ is in $S_{\mathbb R}$, we see that 
$$\rho(S_{\mathbb R})=\rho(S)=\frac14.$$
We further remark that the critical value $\frac14$ coincides with 
\begin{equation}\label{eqn:koebeminusone}
|\mathcal K(-1)|=\min\{|\mathcal K(\zeta)|:\,\zeta\in\mathbb T\},
\end{equation}
where $\mathbb T:=\partial\mathbb D$ is the unit circle, and moreover, $-1$ is the only minimum modulus point for $\mathcal K$ on $\mathbb T$. 

\par Among the many results that highlight the extremal role of the Koebe function $\mathcal K$ in $S$ and/or $S_{\mathbb R}$, the most famous is undoubtedly the (former) {\it Bieberbach conjecture}, now de Branges' theorem, which establishes sharp bounds for the coefficients $a_n$ in \eqref{eqn:taylor}. Namely, it states that every $f\in S$ (and hence every $f\in S_{\mathbb R}$) satisfies $|a_n|\le n$ for  $n=2,3,\dots$, the inequalities being all strict unless $f$ is the Koebe function $\mathcal K$ or one of its rotations. (In the case of $S_{\mathbb R}$, the only nontrivial rotation to be considered is $z\mapsto -\mathcal K(-z)$.) 

\par We mention in passing that the $S_{\mathbb R}$ version of the Bieberbach conjecture was relatively easy to settle; one of the proofs (as outlined in \cite[p.\,269]{Du}) makes use of Suffridge's work on univalent polynomials, a topic to be touched upon below. By contrast, the full version of the conjecture had remained open for almost 70 years, defying numerous attacks, until de Branges finally cracked it by using a highly sophisticated array of techniques (see \cite{deB, FP}). 

\par It is noteworthy that the two extremal problems -- those underlying the Koebe $\frac14$ theorem and the Bieberbach conjecture -- are tightly linked together. In fact, the basic inequality $|a_2|\le 2$, which was discovered by Bieberbach in 1916, both led him to a proof of Theorem \ref{thm:quarter} and provided the basis for his coefficient conjecture. 

\par Now, we are interested in polynomial versions of Theorem \ref{thm:quarter}. The extremal problem that arises, to be described in a moment, will be referred to as the {\it polynomial Koebe problem}. Given a positive integer $N$, let $\mathcal U_N$ (resp., $\mathcal U_{N,\mathbb R}$) denote the set of univalent polynomials $p$ of the form 
\begin{equation}\label{eqn:poly}
p(z)=z+\sum_{n=2}^N a_nz^n
\end{equation}
with complex (resp., real) coefficients; here and below, univalence is only assumed in $\mathbb D$.  Thus, $\mathcal U_N\subset S$ and $\mathcal U_{N,\mathbb R}\subset S_{\mathbb R}$. We then ask: {\it What are the values of the Koebe radii $\rho(\mathcal U_N)$ and $\rho(\mathcal U_{N,\mathbb R})$?} Also, {\it what are the extremal \lq\lq Koebe-type" polynomials that minimize the quantity 
\begin{equation}\label{eqn:dist}
\text{\rm dist}\left(p(\mathbb T),\{0\}\right):=
\min\{|p(\zeta)|:\,\zeta\in\mathbb T\}
\end{equation}
among all $p$ in $\mathcal U_N$ and/or $\mathcal U_{N,\mathbb R}$?}

\par Clearly, the optimal lower bound for \eqref{eqn:dist}, as $p$ ranges over $\mathcal U_N$ or $\mathcal U_{N,\mathbb R}$, coincides with the Koebe radius of the corresponding class. Once again, the case of real coefficients seems to be more tractable, so we restrict most of our attention to $\mathcal U_{N,\mathbb R}$. 

\par Needless to say, the problem is trivial for $N=1$. Indeed, the only element of $\mathcal U_1$ (as well as of $\mathcal U_{1,\mathbb R}$) is the identity function $z$, whence 
$$\rho\left(\mathcal U_1\right)=\rho\left(\mathcal U_{1,\mathbb R}\right)=1.$$
The case $N=2$ is not much harder. This time, the univalent polynomials of the form $z+a_2z^2$ are precisely those with $|a_2|\le\frac12$; the extremal ones have $|a_2|=\frac12$, and a simple calculation shows that 
$$\rho\left(\mathcal U_2\right)=\rho\left(\mathcal U_{2,\mathbb R}\right)=\frac12.$$ Typically enough, passing to higher degrees makes life increasingly painful, and the case of $N=3$ already seems to deserve a serious analysis. To begin with, it is far from trivial to determine the values of $a_2$ and $a_3$ for which the polynomial $z+a_2z^2+a_3z^3$ is univalent in $\mathbb D$. This has been done, however, and we recall the result (or rather its $\mathcal U_{3,\mathbb R}$ version) in Section 4 below. Then we proceed to solve our polynomial Koebe problem for $\mathcal U_{3,\mathbb R}$. 

\par Meanwhile, we pause to speculate on the case of general $N$, trying to guess what the extremal polynomials should be. One natural -- and tantalizingly sexy -- candidate that comes to mind is Suffridge's remarkable family of univalent polynomials, which we discuss in some detail (only to reject it shortly afterwards). 
The Suffridge polynomials are known to approximate and mimic the Koebe function $\mathcal K$ in several ways, so we find it quite surprising that this time they fall short of being extremal. Then we come up with another -- newborn -- collection of polynomials which, we strongly believe, is the right candidate for the job. Finally, the solution we give for $N=3$ serves to corroborate the conjecture, and also allows us to compare the extremal properties of the two competing families of polynomials. 

\section{The Suffridge polynomials -- a rejected candidate}

In \cite{S}, Suffridge introduced an important family of polynomials, which turned out to enjoy a number of elegant extremal properties. Namely, for $N=1,2,\dots$, he defined 
$$q_N(z):=\sum_{k=1}^NA_{k,N}\,z^k,$$
where 
\begin{equation}\label{eqn:defakn}
A_{k,N}:=\frac{N-k+1}N\cdot\frac{\sin\pi k/(N+1)}{\sin\pi/(N+1)},
\end{equation}
and verified that each $q_N$ is univalent in $\mathbb D$. Since $A_{1,N}=1$, we have $q_N\in\mathcal U_{N,\mathbb R}$ for every $N$. We further note that $A_{N,N}=1/N$, 
which already reflects a certain extremal property of $q_N$. (Indeed, the highest coefficient $a_N$ of a univalent polynomial \eqref{eqn:poly} must satisfy $|a_N|\le1/N$. To see why, look at the constant term of the monic polynomial $(Na_N)^{-1}p'(z)$ which has no zeros in $\mathbb D$.) Moreover, Suffridge showed that whenever $p\in\mathcal U_{N,\mathbb R}$ is a polynomial of the form \eqref{eqn:poly} with $|a_N|=1/N$, the remaining coefficients of $p$ are also dominated by those of $q_N$, so that 
\begin{equation}\label{eqn:dom}
|a_k|\le A_{k,N}\quad\text{\rm for}\quad k=2,\dots,N.
\end{equation}
When $N\le4$, it is actually true that every polynomial \eqref{eqn:poly} lying in $\mathcal U_{N,\mathbb R}$ obeys Suffridge's estimates \eqref{eqn:dom}  unrestrictedly, the assumption that $|a_N|=1/N$ being no longer needed. (In the nontrivial cases $N=3$ and $N=4$, this follows from results of \cite{B, CR, K} and of \cite{M}, respectively.) It is also noteworthy that the polynomial 
\begin{equation}\label{eqn:suff3}
q_3(z)=z+\frac{2\sqrt2}3z^2+\frac13z^3
\end{equation}
maximizes $|a_2|$ and $|a_3|$ among all $p$'s of the form $p(z)=z+a_2z+a_3z^3$ in $\mathcal U_3$, not just in $\mathcal U_{3,\mathbb R}$; see \cite{B} or \cite{CR}. 

\par These extremal properties of the Suffridge polynomial $q_N$ seem to indicate that its role in $\mathcal U_{N,\mathbb R}$ is similar to that of the Koebe function $\mathcal K$ in $S$ or $S_{\mathbb R}$, the analogy being especially clear-cut for small degrees. On the other hand, for every fixed $k$, the coefficients $A_{k,N}$ increase to $k$ (i.e., to the $k$th coefficient of the Koebe function) as $N\to\infty$. Consequently, $\lim_{N\to\infty}q_N=\mathcal K$ uniformly on compact subsets of $\mathbb D$. 

\par Now let us try and estimate the distance in \eqref{eqn:dist} for $p=q_N$. The Koebe-type behavior of $q_N$ suggests, in conjunction with \eqref{eqn:koebeminusone}, that we begin by looking at $\left|q_N(-1)\right|$. In fact, Theorem \ref{thm:quarter} tells us that 
\begin{equation}\label{eqn:twosided}
\frac14\le\min\{|q_N(\zeta)|:\zeta\in\mathbb T\}\le\left|q_N(-1)\right|,
\end{equation}
whereas a straightforward computation yields 
$$q_N(-1)=-\frac{N+1}{4N}\left[\cos\frac\pi{2(N+1)}\right]^{-2},$$
so that $\lim_{N\to\infty}q_N(-1)=-\frac14$. We see that the upper bound in \eqref{eqn:twosided} tends to $\frac14$ as $N\to\infty$, meaning that the $q_N$'s are {\it asymptotically sharp} in the polynomial Koebe problem. Are they also sharp for each individual $N$? 

\par This last question was explicitly raised in \cite{Dim}, where polynomial analogues of Theorem \ref{thm:quarter} were also touched upon, and this has largely spurred our interest in the problem. While the above discussion seems to provide evidence in favor of a \lq\lq yes" answer, we now disprove the conjecture (at least in the $\mathcal U_{N,\mathbb R}$ setting) by showing that the actual answer is a resounding \lq\lq no," already for $N=3$. As a matter of fact, the Suffridge polynomial \eqref{eqn:suff3} fails to be extremal for the Koebe problem in $\mathcal U_{3,\mathbb R}$ since it loses the game to 
\begin{equation}\label{eqn:ourp3}
p_3(z):=z+\frac2{\sqrt5}z^2+\frac12\left(1-\frac1{\sqrt5}\right)z^3,
\end{equation}
another remarkable polynomial from $\mathcal U_{3,\mathbb R}$, which turns out to be unbeatable. Specifically, the corresponding values of the distance in \eqref{eqn:dist} happen to be 
\begin{equation}\label{eqn:distourp3}
\min_{\zeta\in\mathbb T}\left|p_3(\zeta)\right|=\left|p_3(-1)\right|
=\frac{3-\sqrt5}2=0.3819\dots
\end{equation}
and 
\begin{equation}\label{eqn:distsuff3}
\min_{\zeta\in\mathbb T}\left|q_3(\zeta)\right|
=\left|q_3\left(-\frac{2\sqrt2}3\pm\frac i3\right)\right|=\frac2{3\sqrt3}
=0.3849\dots,
\end{equation}
so $p_3$ does indeed slightly better. 

\par The facts just mentioned (i.e., the univalence of $p_3$, its extremality in the Koebe problem for $\mathcal U_{3,\mathbb R}$, and hence also its supremacy over $q_3$) will be verified in Sections 4 and 5 below. The calculations leading to \eqref{eqn:distourp3} and \eqref{eqn:distsuff3} will be provided there as well. But first we have to place the polynomial $p_3$ where it belongs. Namely, it should be viewed as a member of a certain lordly family, $\{p_N\}$, which we now describe. 

\section{A new, more promising, family of polynomials\\ 
and the conjectured solution}

Recall, to begin with, that the Chebyshev polynomials of the second kind, $U_n$, are defined for $n=0,1,2,\dots$ by the identity 
$$U_n(\cos\theta)=\frac{\sin(n+1)\theta}{\sin\theta},\qquad\theta\in(-\pi,\pi].$$
Thus, 
$$U_0(t)=1,\quad U_1(t)=2t,\quad U_2(t)=4t^2-1,\quad U_3(t)=8t^3-4t,$$
and so forth. Next, for a positive integer $N$, we put 
$$c_N:=\cos\frac\pi{N+2}$$ 
and consider the numbers
\begin{equation}\label{eqn:defbkn}
B_{k,N}:=\frac{U'_{N-k+1}(c_N)}{U'_N(c_N)}\cdot U_{k-1}(c_N)
\end{equation}
with $k=1,\dots,N$. Finally, we define 
$$p_N(z)=\sum_{k=1}^NB_{k,N}\,z^k\qquad (N=1,2,\dots).$$
It should be noted that $B_{1,N}=1$ for each $N$. Also, rewriting the expression \eqref{eqn:defakn} for the Suffridge coefficients as 
$$A_{k,N}:=\frac{N-k+1}N\cdot U_{k-1}(c_{N-1}),$$
one might observe a certain -- perhaps remote -- kinship between the $A_{k,N}$ and the $B_{k,N}$, or equivalently, between the two families of polynomials. The new formulas \eqref{eqn:defbkn} look somewhat more bizarre, if not a bit scary, but there are reasons for them being what they are. 

\par In fact, the polynomials $p_N$ arose quite recently (see \cite{DSS}) in connection with another extremal problem, which is fairly close in spirit to the current one. The problem was: Given $N\in\mathbb N$, maximize the quantity 
\begin{equation}\label{eqn:mup}
\mu(p):=\min\left\{\text{\rm Re}\,p(\zeta):\,\zeta\in\mathbb T,\,
\text{\rm Im}\,p(\zeta)=0\right\}
\end{equation}
over all polynomials $p$ of the form \eqref{eqn:poly} with real coefficients (but without assuming univalence). It was then shown in \cite{DSS} that the unique maximizing polynomial is precisely $p_N$, so that the best upper bound for $\mu(p)$ is $\mu(p_N)$, which in turn equals $-1/(4c_N^2)$. 

\par The first two polynomials in the $p_N$ family are 
$$p_1(z)\left(=q_1(z)\right)=z$$
and 
$$p_2(z)\left(=q_2(z)\right)=z+\frac12z^2,$$
both being obviously univalent in $\mathbb D$. The next one, $p_3$, is our old friend \eqref{eqn:ourp3} which is again univalent in $\mathbb D$, as we shall see in Section 4 below. Then comes 
$$p_4(z)=z+\frac76z^2+\frac23z^3+\frac16z^4,$$
a polynomial whose univalence has also been established; a nice proof can be found in \cite{Dil}. In fact, $p_4$ is even known to be starlike (meaning that it maps $\mathbb D$ conformally onto a starlike domain), since it meets the starlikeness criterion given in \cite[pp.\,515--516]{GH}. We have been able to verify univalence for $p_5$ and $p_6$ as well, but the case of bigger $N$'s remains open. We do believe that $p_N$ is actually univalent in $\mathbb D$, and hence $p_N\in\mathcal U_{N,\mathbb R}$, for all $N$. Numerical simulations reinforce this belief substantially. 

\par We further conjecture that the $p_N$'s are extremal in the polynomial Koebe problem, so that for every fixed $N$ in $\mathbb N$, $p_N$ minimizes the distance in \eqref{eqn:dist} among all $p\in\mathcal U_{N,\mathbb R}$. The Koebe radius $\rho(\mathcal U_{N,\mathbb R})$ must then agree with $\min_{\zeta\in\mathbb T}|p_N(\zeta)|$, and it is very likely that this last quantity always equals $|p_N(-1)|$, which in turn simplifies to $1/(4c_N^2)$. We go on to claim that the same result should hold in the case of complex coefficients, so that $\mathcal U_N$ has presumably the same Koebe radius and the same extremal polynomials as $\mathcal U_{N,\mathbb R}$. Thus, in particular, it is conjectured that 
\begin{equation}\label{eqn:conjrad}
\rho\left(\mathcal U_N\right)=\rho\left(\mathcal U_{N,\mathbb R}\right)
=\frac1{4\cos^2(\pi/(N+2))},\qquad N\in\mathbb N.
\end{equation}

\par There are several sources for our certainty, beyond a reasonable doubt, that the conjectured solution is correct. These include the appearance of the $p_N$ polynomials in the cognate extremal problem involving \eqref{eqn:mup}, as mentioned above, plus the analysis of the case $N=3$ (which is the bifurcation point between the $p_N$'s and $q_N$'s) to be carried out below, plus the numerical experiments we performed when playing around with polynomials of higher degrees. 

\section{Polynomials of degree 3: preliminaries} 

Which functions $p$ of the form 
\begin{equation}\label{eqn:cubicpoly}
p(z)=z+a_2z^2+a_3z^3 
\end{equation}
are univalent in $\mathbb D$? The answer was first obtained in \cite{K} and then rediscovered, via different approaches, in \cite{B} and \cite{CR}. To state it (which we only do for the case where $a_2$ and $a_3$ are real), we begin by describing the boundary $\Gamma$ of the univalence region in the $(a_2,a_3)$ plane. 

\par The portion $\Gamma_+$ of $\Gamma$ that lies in the half-plane $\{a_2\ge0\}$ can be written as $\Gamma_+=\gamma_1\cup\gamma_2\cup\gamma_3$, where 
\begin{itemize}
\item $\gamma_1$ is the segment of the line $2a_2-3a_3=1$ with endpoints $\left(0,-\frac13\right)$ and $\left(\frac45,\frac15\right)$; 
\item $\gamma_2$ is the (shorter) arc of the ellipse $a_2^2+4(a_3-\frac12)^2=1$ with endpoints $\left(\frac45,\frac15\right)$ and $\left(\frac{2\sqrt2}3,\frac13\right)$; 
\item $\gamma_3$ is the segment of the line $a_3=\frac13$ with endpoints $\left(\frac{2\sqrt2}3,\frac13\right)$ and $\left(0,\frac13\right)$.
\end{itemize}
(The two line segments and the arc are assumed to be closed.) We then define 
$$\Gamma_-:=\left\{(a_2,a_3)\in\mathbb R^2:\,(-a_2,a_3)\in\Gamma_+\right\}$$ 
and $\Gamma:=\Gamma_+\cup\Gamma_-$. Thus, $\Gamma$ is a simple closed curve which is symmetric with respect to the $a_3$ axis. Finally, we write $\Omega$ for the bounded connected component of $\mathbb C\setminus\Gamma$ and put $V:=\Omega\cup\Gamma$.

\par The required result from (any of) \cite{B,CR,K} can now be stated as follows. 

\begin{lem}\label{lem:univcrit} For $(a_2,a_3)\in\mathbb R^2$, the polynomial \eqref{eqn:cubicpoly} is univalent in $\mathbb D$ if and only if $(a_2,a_3)$ belongs to $V$.
\end{lem}

\par Obviously enough, the univalence region $V$ is also symmetric with respect to the $a_3$ axis (as is $\Gamma$). This is due to the fact that the polynomial \eqref{eqn:cubicpoly} and its reflection 
$$p^*(z):=-p(-z)=z-a_2z^2+a_3z^3$$
are, or are not, univalent simultaneously. 

\par The Suffridge polynomial \eqref{eqn:suff3} corresponds to the vertex $\left(\frac{2\sqrt2}3,\frac13\right)=\gamma_2\cap\gamma_3$ which sticks out in both coordinate directions and has a special, remarkably extreme position in $V$ (along with the symmetric point $\left(-\frac{2\sqrt2}3,\frac13\right)$ that represents $q^*_3$). Thus, when faced with an extremal problem for $\mathcal U_{3,\mathbb R}$, one is indeed tempted to contemplate $q_3$ and $q^*_3$ as the most likely extremizers. In our case, however, the actual winners turn out to be $p_3$ (as defined by \eqref{eqn:ourp3} above) and $p^*_3$, a fact we shall soon verify. The corresponding points in the coefficient plane are 
\begin{equation}\label{eqn:ourpoints}
\left(\pm\frac2{\sqrt5},\frac12\left(1-\frac1{\sqrt5}\right)\right);
\end{equation}
both belong to $\Gamma$ (in fact, the one with the $+$ sign lies on the arc $\gamma_2$), so univalence is ensured by Lemma \ref{lem:univcrit}. The points \eqref{eqn:ourpoints} do not appear to enjoy a particularly privileged position, though, so the extremal nature of $p_3$ and $p^*_3$ can scarcely be viewed as predictable. 

\par Another preliminary question we need to discuss is this: Given a polynomial
\eqref{eqn:cubicpoly} with $a_2,a_3\in\mathbb R$, where does it attain its minimum modulus value on the unit circle $\mathbb T$? 

\par We are only interested in the case $a_3\ne0$. For a point $z=x+iy\in\mathbb T$, a straightforward calculation yields
\begin{equation}\label{eqn:modpsq}
\begin{aligned}
|p(z)|^2
&=(1+a_2z+a_3z^2)(1+a_2\overline z+a_3\overline z^2)\\
&=1+a_2^2+a_3^2-2a_3+2a_2(1+a_3)x+4a_3x^2=:\Phi(x).
\end{aligned}
\end{equation}
Originally, $x:=\text{\rm Re}\,z$ runs through the interval $[-1,1]$, but we extend the quadratic polynomial $\Phi(x)$ to all $x\in\mathbb R$. Then 
$$\Phi'(x)=2a_2(1+a_3)+8a_3x,$$
and the only zero of this derivative is 
\begin{equation}\label{eqn:xnought}
x_0=-\frac{a_2(1+a_3)}{4a_3}.
\end{equation}
The function $\Phi(x)$ therefore attains its minimum (if $a_3>0$) or maximum (if $a_3<0$) at $x_0$, its value at the critical point being 
\begin{equation}\label{eqn:phicrit}
\Phi(x_0)=(1-a_3)^2\left(1-\frac{a_2^2}{4a_3}\right),
\end{equation}
as verified by direct computation. In particular, this last quantity will be nonnegative whenever $x_0$ happens to be in $[-1,1]$; to see why, recall \eqref{eqn:modpsq}. 

\par In terms of $p$, two types of behavior may occur. To distinguish between them, we now introduce the appropriate terminology. 

\begin{defn} A polynomial $p$ is said to be {\it of type I} if 
\begin{equation}\label{eqn:typeone}
\min\{|p(z)|:\,z\in\mathbb T\}=\min\{|p(-1)|,|p(1)|\}.
\end{equation}
Otherwise we say that $p$ is {\it of type II}.
\end{defn}

\par It should be noted that, for a polynomial $p$ with nonnegative coefficients, \eqref{eqn:typeone} simplifies to
\begin{equation}\label{eqn:typeonesim}
\min\{|p(z)|:\,z\in\mathbb T\}=|p(-1)|.
\end{equation}
Thus, polynomials of type I are essentially those that mimic the Koebe function $\mathcal K$ by sharing its property \eqref{eqn:koebeminusone}. 

\par The above discussion leads us to the following conclusion. 

\begin{lem}\label{lem:typetwo} Let $p$ be a polynomial of the form \eqref{eqn:cubicpoly} with real coefficients and with $a_3\ne0$. In order that $p$ be of type II, it is necessary and sufficient that $a_3>0$ and $-1<x_0<1$, where $x_0$ is defined by \eqref{eqn:xnought}. In this case, 
\begin{equation}\label{eqn:mintypetwo}
\min\{|p(z)|:z\in\mathbb T\}=
|p(x_0\pm iy_0)|=|1-a_3|\left(1-\frac{a_2^2}{4a_3}\right)^{1/2},
\end{equation}
where $y_0:=\sqrt{1-x_0^2}$. Moreover, $x_0\pm iy_0$ are then the only points of $\mathbb T$ where the minimum in question is attained.
\end{lem}

\par This result allows us to compute the quantity \eqref{eqn:dist} for $p=p_3$ and $p=q_3$, thus verifying the announced formulas \eqref{eqn:distourp3} and 
\eqref{eqn:distsuff3}. 

\par The polynomial $p_3$ has $a_2=\frac2{\sqrt5}$ and $a_3=\frac12\left(1-\frac1{\sqrt5}\right)$, and plugging this into \eqref{eqn:xnought} gives 
$$x_0=-\frac14-\frac7{20}\sqrt5=-1.0326\dots.$$ 
It now follows from Lemma \ref{lem:typetwo} that $p_3$ is of type I, and so 
\begin{equation}\label{eqn:prepdistourp3}
\min\{|p_3(z)|:\,z\in\mathbb T\}=|p_3(-1)|.
\end{equation}
The right-hand side of \eqref{eqn:prepdistourp3} reduces to $(3-\sqrt5)/2$, and we arrive at \eqref{eqn:distourp3}. 

\par As to the Suffridge polynomial $q_3$, this time we have 
$$x_0=-\frac{2\sqrt2}3(=-a_2)=-0.9428\dots,$$
so Lemma \ref{lem:typetwo} tells us that $q_3$ is of type II. The corresponding $y_0$ equals $\frac13(=a_3)$, and substituting the appropriate values into \eqref{eqn:mintypetwo} yields \eqref{eqn:distsuff3}. 

\section{Polynomials of degree 3: solution} 

For a polynomial $p$, we put 
$$m(p):=\min\{|p(\zeta)|:\,\zeta\in\mathbb T\}.$$ 
To solve the Koebe problem for $\mathcal U_{3,\mathbb R}$, we need to minimize the functional $m(p)$ over all $p\in\mathcal U_{3,\mathbb R}$. This is done in Theorem \ref{thm:classification} below, where the minimizing polynomials are exhibited; as promised, these are shown to be $p_3$ and $p_3^*$. 

\par A couple of conventions will be made. First, if $X$ is a class of polynomials and $F\in X$, we say that $F$ {\it is extremal for $X$} to mean that 
$$\inf\{m(p):\,p\in X\}=m(F).$$
Secondly, every polynomial $p$ in $\mathcal U_{3,\mathbb R}$ will be identified, via \eqref{eqn:cubicpoly}, with the corresponding point $(a_2,a_3)$ in the plane (or rather in the univalence region $V$ coming from Lemma \ref{lem:univcrit}); we shall occasionally write $p=(a_2,a_3)$ to make this explicit. Also, given a set $M\subset V(\subset\mathbb R^2)$, we may now use the notation $p\in M$ without any risk of confusion. 

\begin{thm}\label{thm:classification} The only extremal polynomials for the class $\mathcal U_{3,\mathbb R}$ are $p_3$, as defined by \eqref{eqn:ourp3}, and its reflection $p_3^*$. 
\end{thm}

\par As a consequence, we see that the Koebe radius for $\mathcal U_{3,\mathbb R}$ equals $m(p_3)$, which agrees with 
$$\frac1{4\cos^2\frac\pi5}=\frac{3-\sqrt5}2,$$ 
the conjectured (and now established) value of $\rho(\mathcal U_{3,\mathbb R})$ from \eqref{eqn:conjrad}. It only remains to prove Theorem \ref{thm:classification}. 

\begin{proof} The extremal polynomials must live on the boundary, $\Gamma$, of the univalence region $V$. By symmetry, it suffices to consider
$$\Gamma_+=\{(a_2,a_3)\in\Gamma:\,a_2\ge0\},$$
which in turn decomposes as $\gamma_1\cup\gamma_2\cup\gamma_3$; see the preceding section for definitions. 
\par We begin by looking at the values of $m(p)$ when $p\in\gamma_1$. It is easy to check that, whenever $p=(a_2,a_3)$ is a point of $\gamma_1$ with $a_3>0$, we have 
$$\frac{a_2(1+a_3)}{4a_3}\ge1.$$
Equivalently, the number $x_0$ given by \eqref{eqn:xnought} is in $(-\infty,-1]$ for any such point, and we deduce from Lemma \ref{lem:typetwo} that the polynomials belonging to $\gamma_1$ are all of type I. It is also clear that every polynomial $p\in\gamma_1$ satisfies 
\begin{equation}\label{eqn:minusplusone}
|p(-1)|=1-a_2+a_3\le1+a_2+a_3=|p(1)|.
\end{equation}
Furthermore, on $\gamma_1$ we have $a_3=\frac13(2a_2-1)$, whence 
\begin{equation}\label{eqn:lineeq}
1-a_2+a_3=\frac13(2-a_2),
\end{equation}
the permissible values of $a_2$ being those in $[0,\frac45]$. We now combine \eqref{eqn:typeone}, \eqref{eqn:minusplusone}, and \eqref{eqn:lineeq} to get 
$$m(p)=|p(-1)|=\frac13(2-a_2)$$
for each polynomial $p=(a_2,a_3)\in\gamma_1$. Consequently, 
\begin{equation}\label{eqn:infgammaone}
\inf\{m(p):\,p\in\gamma_1\}=\frac13\left(2-\frac45\right)=\frac25=0.4.
\end{equation}

\par Next, we turn to the case where $p=(a_2,a_3)\in\gamma_2$. We know that $\gamma_2$ contains polynomials of both types, and we write $\sigma_{\rm I}$ (resp., $\sigma_{\rm II}$) for the set of all polynomials of type I (resp., II) that are in $\gamma_2$. In fact, $\sigma_{\rm I}$ and $\sigma_{\rm II}$ are disjoint subarcs of $\gamma_2$ whose common endpoint $\widetilde p=(\widetilde a_2,\widetilde a_3)\in\gamma_2$ is determined by the relation 
$$\widetilde a_2=\frac{4\widetilde a_3}{1+\widetilde a_3}.$$
(Thus, plugging $a_2=\widetilde a_2$ and $a_3=\widetilde a_3$ into \eqref{eqn:xnought} yields $x_0=-1$. A bit of inspection shows that $\gamma_2$ contains exactly one point with this property. We note that the number $\widetilde a_3$ satisfies $\frac15<\widetilde a_3<\frac13$ and coincides with the unique positive root of the equation $t^3+t^2+3t-1=0$.) Precisely speaking, $\sigma_{\rm I}$ is the closed subarc of $\gamma_2$ with endpoints $\left(\frac45,\frac15\right)$ and $(\widetilde a_2,\widetilde a_3)$, while 
$\sigma_{\rm II}=\gamma_2\setminus\sigma_{\rm I}$ is the complementary (half-open) subarc with endpoints $(\widetilde a_2,\widetilde a_3)$ and $\left(\frac{2\sqrt2}3,\frac13\right)$. 

\par Rewriting the equation of the ellipse (of which $\gamma_2$ forms part) as 
\begin{equation}\label{eqn:ellipse}
a_2^2=4a_3(1-a_3),
\end{equation}
we may parametrize $\gamma_2$ in the form 
$$\gamma_2=\left\{\left(2\sqrt{t(1-t)},\,t\right):\,
t\in\left[\frac15,\frac13\right]\right\},$$
where $\sigma_{\rm I}$ and $\sigma_{\rm II}$ correspond to the 
parameter ranges $\left[\frac15,\widetilde a_3\right]=:J_{\rm I}$ and $\left(\widetilde a_3,\frac13\right]=:J_{\rm II}$. 

\par Accordingly, the quantity $m(p)$ with $p\in\gamma_2$ admits a fairly simple expression in terms of the $a_3$ coordinate. Namely, for $p=(a_2,a_3)\in\sigma_{\rm I}$, we combine \eqref{eqn:typeonesim} with \eqref{eqn:ellipse} to obtain
\begin{equation}\label{eqn:sigone}
m(p)=|p(-1)|=1-a_2+a_3=1-2\sqrt{a_3(1-a_3)}+a_3,
\end{equation}
whereas for $p=(a_2,a_3)\in\sigma_{\rm II}$ we invoke \eqref{eqn:mintypetwo} in conjunction with \eqref{eqn:ellipse} to get 
\begin{equation}\label{eqn:sigtwo}
m(p)=(1-a_3)\left(1-\frac{a_2^2}{4a_3}\right)^{1/2}=\sqrt{a_3}(1-a_3).
\end{equation}
Differentiating, we find that the function 
$$\varphi(t):=1-2\sqrt{t(1-t)}+t,\qquad t\in J_{\rm I},$$
has a minimum at 
$$t_*:=\frac12\left(1-\frac1{\sqrt5}\right)$$
(which is an interior point of $J_{\rm I}$) and  
$$\min\left\{\varphi(t):\,t\in J_{\rm I}\right\}=\varphi(t_*)=\frac{3-\sqrt5}2.$$ 
Moreover, $t_*$ is the only point in $J_{\rm I}$ where the minimum in question is attained. We further observe that the function 
$$\psi(t):=\sqrt t(1-t),\qquad t\in J_{\rm II},$$
is increasing on its domain $J_{\rm II}(\subset(0,\frac13])$; this is again verified by differentiation. 

\par From \eqref{eqn:sigone} and \eqref{eqn:sigtwo} we know that $m(p)$ equals $\varphi(a_3)$ when $p\in\sigma_{\rm I}$, and $\psi(a_3)$ when $p\in\sigma_{\rm II}$. The critical value $t_*$ of the $a_3$ variable corresponds to the point $(2/\sqrt5,t_*)\in\sigma_{\rm I}$, which represents the polynomial $p_3$. Because $m(p)$ is also continuous at $\widetilde p$ (as it is everywhere else), the above facts about the $\varphi$ and $\psi$ functions allow us to conclude that 
\begin{equation}\label{eqn:infgammatwo}
\inf\left\{m(p):\,p\in\gamma_2\right\}=m(p_3)=\frac{3-\sqrt5}2
\end{equation}
and that $p_3$ is the only extremal polynomial for $\gamma_2$. 

\par It remains to consider the case where $p=(a_2,a_3)\in\gamma_3$. Since $0\le a_2\le\frac{2\sqrt2}3$ and $a_3=\frac13$, the formula \eqref{eqn:xnought} yields $x_0=-a_2$, whence in particular 
$$-1<-\frac{2\sqrt2}3\le x_0\le0.$$ 
It now follows from Lemma \ref{lem:typetwo} that $p$ is of type II, and moreover, 
$$m(p)=\frac13\left(4-3a_2^2\right)^{1/2}.$$
Clearly, this is minimized by assigning the largest admissible value, $\frac{2\sqrt2}3$, to the $a_2$ variable. In other words, the only extremal polynomial for $\gamma_3$ is $q_3=\left(\frac{2\sqrt2}3,\frac13\right)$, and so 
\begin{equation}\label{eqn:infgammathree}
\inf\{m(p):\,p\in\gamma_3\}=m(q_3)=\frac2{3\sqrt3}.
\end{equation}

\par Finally, a quick glance at \eqref{eqn:infgammaone}, \eqref{eqn:infgammatwo} and \eqref{eqn:infgammathree} reveals that the smallest of the three infima is that over $\gamma_2$. This implies the extremality (and uniqueness) of $p_3$ among all polynomials $p=(a_2,a_3)\in\mathcal U_{3,\mathbb R}$ with $a_2\ge0$. By symmetry, a similar role is played by $p^*_3$ among the $\mathcal U_{3,\mathbb R}$ polynomials with $a_2<0$. The proof is complete. 
\end{proof}

\end{document}